\documentclass[12pt]{amsart}
\usepackage{amsmath,amssymb,amsthm}
\usepackage[]{latexsym}
\newtheorem{theorem}{Theorem}[section]
\newtheorem{definition}[theorem]{Definition}

\newtheorem{proposition}[theorem]{Proposition}

\newtheorem{lemma}[theorem]{Lemma}

\newcommand{\gt}{{\tilde {\gamma }}}
\newcommand{\et}{{\tilde {\varepsilon }}}
\newcommand{\gtau}{{\gamma _{\tau }}}
\newcommand{\etau}{{\varepsilon _{\tau}}}
\newcommand{\om}{{\omega }}

\newcommand{\cR}{{\mathcal R}}

\newcommand{\cD}{{\mathcal D}}

\newcommand{\cH}{{\mathcal H}}

\newcommand{\cL}{{\mathcal L}}

\newcommand{\cO}{{\mathcal O}}

\newcommand{\Pol}{{\mathrm {Pol }}}

\newcommand{\sign}{{\mathrm {sign }}}

\newcommand{\Div}{{\mathrm {Div }}}

\newcommand{\Gal}{{\mathrm {Gal }}}
\newcommand{\disc}{{\mathrm {disc }}}

\newcommand{\Pic}{\mathrm{Pic }}

\newcommand{\Norm}{\mathrm{N }}
\newcommand{\M}{\mathrm{M }}
\newcommand{\Aut}{\mathrm{Aut }}
\newcommand{\PGL}{{\mathrm \mathbb{PGL}}}
\newcommand{\GL}{{\mathrm \mathbb{GL}}}
\newcommand{\End}{{\mathrm{End }}}

\newcommand{\NS}{{\mathrm{NS }}}
\newcommand{\n}{{\mathrm{n}}}

\newcommand{\rank}{{\mathrm {rank }}}
\newcommand{\Fr}{{\mathrm {Frob }}}
\newcommand{\tr}{{\mathrm {tr }}}

\newcommand{\Res}{{\mathrm Res}}

\newcommand{\cond}{{\mathrm {cond }}}

\newcommand{\Z}{{\mathbb Z}}
\newcommand{\Q}{{\mathbb Q}}
\newcommand{\C}{{\mathbb C}}
\newcommand{\R}{{\mathbb R}}
\newcommand{\F}{{\mathbb F}}

\newcommand{\PP}{{\mathbb P}}

\newcommand{\ra}{\rightarrow}
\newcommand{\lra}{\longrightarrow}
\newcommand{\qbar}{{\overline {\Q }}}

\newcommand{\Ol}{\mathcal{O}_\ell^*}
\newcommand{\Olt}{\widetilde{\mathcal{O}}_\ell^*}
\newcommand{\la}{\lambda}
\newcommand{\U}{\mathcal{U}_\lambda}
\newcommand{\V}{\mathcal{V}_\lambda}

\newcommand{\re}{\rho_\lambda}
\newcommand{\rer}{\bar{\rho}_\lambda}
\newcommand{\tres}{\breve{3}}
\newcommand{\Image}{{\rm Image}}

\include{thebibliography}

\begin{document}

\title{The
arithmetic of QM-abelian surfaces \\ through their Galois
representations}

\author{Luis V. Dieulefait, Victor Rotger}

\footnote{Supported by
 European Research Networks at the Universit\'{e} de Paris 13 and the
 Institut de Mat\'{e}matiques de Jussieu
and by a MECD postdoctoral grant at the Centre de Recerca
Matem\'{a}tica from Ministerio de Educaci\'{o}n y Cultura;
Partially supported by DGICYT Grant BFM2003-06768-C02-02. }

\address{Centre de Recerca Matem\'{a}tica, Apartat 50, E-08193 Bellaterra, Spain;
 Escola Universit\`aria Polit\`ecnica de
Vilanova i la Geltr\'u, Av. V\'{\i}ctor Balaguer s/n, E-08800 Vilanova
i la Geltr\'u, Spain}

\email{LDieulefait@crm.es, vrotger@mat.upc.es}



\begin{abstract}

This note provides an in\-sight to the diophan\-tine proper\-ties
of abelian surfaces with quaternionic multiplication over number
fields. We study the fields of definition of the endomorphisms on
these abelian varieties and the images of the Galois
representations on their Tate modules. We illustrate our results
with an explicit example.

\end{abstract}

\maketitle
\section{Abelian surfaces with quaternionic multiplication}
\markboth{Luis V. Dieulefait, Victor Rotger}{Arithmetic of abelian
surfaces with quaternionic multiplication}

Fix $\qbar $ an algebraic closure of the field $\Q $ of rational
numbers and let $K\subset \qbar $ be a number field. Let $A$ be an
abelian surface defined over $K$. Due to Albert's classification
of involuting division algebras (cf. \cite{Mu}), there is a
limited number of possible structures for the algebra of
endomorphisms $\End _{\qbar }(A)\otimes \Q $ of $A$.

We focus our attention on the quaternionic case. While the
existing literature concerning the theme mainly restrict to
abelian surfaces with multiplication by a maximal order in a
division quaternion algebra, in this note we consider quaternionic
multiplication in a wider sense: we shall assume that $\End
_{\qbar }(A)\otimes \Q $ is an arbitrary indefinite quaternion
algebra $B$ over $\Q $, including the split case $B=\M _2(\Q )$,
which is recurringly encountered in the modular setting. This
allows the abelian surface to be isogenous to the product of two
isogenous elliptic curves without CM.

Moreover, we will let $\cO =\End _{\qbar }(A)$ be an arbitrary
order in $B$, although our main results restrict to so called
hereditary orders. We need to be careful on the exact order $\End
_{\qbar }(A)\subset B$ of endomorphisms of $A$ since we are
interested on properties of $A$ that heavily depend on its
isomorphism class and badly behave up to isogeny.

Let then $B = (\frac {a, b}{\Q })= \Q + \Q i + \Q j + \Q i j$, $i
j = -j i$, $i^2=a$, $j^2=b$ with $a, b\in \Q ^*$, be a quaternion
algebra and let $\tr :B\ra \Q $ and $\n :B\ra \Q $, denote the
reduced trace and norm, respectively. The algebra $B$ is said to
be indefinite if the archimedean place of $\Q $ is unramified:
$B\otimes \R \simeq \M _2(\R )$. Equivalently, $B$ is indefinite
if either $a>0$ or $b>0$.

An order $\cO $ in $B$ is a subring of rank $4$ over $\Z $. It is
called a {\em maximal} order if $\cO $ is not properly contained
in any other and it is an {\em Eichler} order if $\cO = \cO _1\cap
\cO _2$ is the intersection of two maximal orders $\cO _1$, $\cO
_2$ in $B$. An order $\cO $ is {\em hereditary} if all its
one-sided modules are projective.

The (reduced) discriminant of an order is $\disc (\cO ) = |\det
(\tr (\beta _i \bar {\beta _j}))|^{1/2}$ for any $\Z $-basis $\{
\beta _1, \beta _2, \beta _3, \beta _4\} $ of $\cO $. The
discriminant of a maximal order is square-free and, since only
depends on $B$, it is simply denoted $\disc (B)$. We have that
$p\mid \disc (B)$ if, and only if, $B_p=B\otimes \Q _p$ is a
division algebra over $\Q _p$. If $\cO $ is an Eichler order, then
$\disc (\cO )= \disc (B)\cdot N$ for some $N>0$ coprime to $\disc
(B)$; $N$ is called the level of $\cO $. Hereditary orders are
exactly the orders in $B$ of square-free discriminant. It can be
shown that an order $\cO \subset B$ is hereditary if and only if
it is an Eichler order of square-free level.

Under the indefiniteness assumption, all one-sided ideals of an
hereditary order are principal. Moreover, two hereditary orders
$\cO $, $\cO '$ in $B$ are isomorphic if and only if $\disc (\cO
)=\disc (\cO ')$. We refer the reader to \cite{AlBa}, \cite{Joh},
\cite{Re} and \cite{Vi} for more details on quaternion algebras
and orders.

\begin{definition}

Let $\cO $ be an order in an indefinite quaternion algebra $B$
over $\Q $. An abelian surface has quaternionic multiplication by
$\cO $ if there is an isomorphism $\iota : \cO \stackrel {\sim
}{\rightarrow } \End _{\qbar }(A)$. A field of definition for the
pair $(A, \iota )$ is an extension $L/K$ such that $\iota :\cO
\stackrel {\sim }{\rightarrow }\End _L(A)$.

\end{definition}

It is one of the aims of this paper to study

\begin{enumerate}
\item
The field extension $L/K$ given by the field of definition $L$ of
the quaternionic multiplication on an abelian surface $A/K$.

\item
The filtration of intermediate endomorphism algebras $\End
_E(A)\otimes \Q \subseteq B$ for $K\subseteq E\subseteq L$.
\end{enumerate}

The first question was studied in greater generality by A.
Silverberg in (\cite{Si}) and Ribet in \cite{Ri2}. When
particularized to our situation we obtain the first interesting
result on this direction.

\begin{proposition}\cite{Si}, \cite{Ri2}
\label{teo:SiRi}

Let $A/K$ be an abelian variety over a number field $K$ and let
$\cO \subseteq \End _{\qbar }(A)$ be a subring of endomorphisms of
$A$. Then there is a {\em unique minimal extension} $L/K$ such
that $\cO \subseteq \End _L(A)$.

The extension $L/K$ is normal and non-ramified at the prime ideals
of $K$ of good or semistable reduction of $A$.

\end{proposition}

With respect to proposition \ref{teo:SiRi}, let us remark that, if
$\End _{\qbar }(A)$ is an order in a division quaternion algebra,
then $A/K$ has potential good reduction and therefore no places of
$K$ of bad reduction of $A$ are semistable. This is a consequence
of Grothendieck's potential good reduction theorem.

Further, Silverberg gave an upper bound for the degree $[L:K]$ in
terms of certain combinatorial numbers. In the particular case of
abelian surfaces with quaternionic multiplication, prop. 4.3 of
\cite{Si} predicts that $[L:K]\leq 48$. As our results will show,
these bounds are not sharp (see proposition \ref{ChFr} for
arbitrary orders $\cO $ and theorem \ref{LL} for hereditary
orders).

Concerning the second question, the non-trivial sub-algebras of
$B$ over $\Q $ are the quadratic fields $\Q (\sqrt {d})$ for $d\in
\Z $ such that any prime number $p\mid D$ does not split in $\Q
(\sqrt {d})$. Although there are infinitely many choices of them,
main theorem \ref{LL} shows that, under the assumption on $\cO $
to be hereditary, there are very restrictive conditions for $\Q
(\sqrt{d})$ to be realized as the algebra of endomorphisms  of $A$
over $K$.

Below, for any positive integer $N$, we write $\M _0(N)=\{ \begin{pmatrix} a&b\\
c N & d\end{pmatrix}, a, b, c, d\in \Z \} $ for the matrix Eichler
order of level $N$. When particularized to the Jacobian variety of
a curve $C/K$ of genus $2$, theorem \ref{LL} asserts the
following.

\begin{theorem}
\label{L}
Let $C/K$ be a curve of genus $2$ defined over a number field $K$
and let $J(C)$ be its Jacobian variety.

{\bfseries I. }[Simple case] Assume that $J(C)$ is absolutely
simple and that $\End _{\qbar }(J(C))=\cO $ is an hereditary order
of discriminant $D=\disc (\cO )$ in a quaternion algebra $B$.

Let $L/K$ be the minimal extension of $K$ over which all
endomorphisms of $J(C)$ are defined. Then

\begin{enumerate}

\item
$L/K$ is an abelian extension with $G = \Gal (L/K)\simeq (1)$,
$C_2$ or $D_2 = C_2\times C_2$, where $C_2$ denotes the cyclic
group of order two.

\item
If $B \not \simeq (\dfrac {-D, m}{\Q })$ for any $m\mid D$, then
$L/K$ is at most a quadratic extension of $K$. In this case, $\End
_K(A)\simeq \Q (\sqrt {-D})$.

\item
If $B = (\dfrac {-D, m}{\Q })$ for some $m\mid D$, then $\End
_K(A)$ is isomorphic to either $\cO $, an order in $\Q (\sqrt
{-D})$, $\Q (\sqrt {m} )$ or $\Q (\sqrt {D/m})$, or $\Z $. In each
case, we respectively have $\Gal (L/K)\simeq (1)$, $C_2$ and
$D_2$.

\end{enumerate}

{\bfseries II. }[Split case] Assume that there is an isomorphism
$\psi _{/\qbar }:J(C)\stackrel {\sim }{\ra }E_1\times E_2$ of
$J(C)$ onto the product of two isogenous elliptic curves $E_1$,
$E_2$ without CM over $\qbar $. Let $\varphi :E_1\ra E_2$ be an
isogeny of minimal degree between them and assume that $N=\deg
(\varphi )$ is square-free. Let $L=\Q (\varphi , \psi )$ be the
compositum of the minimal fields of definition of $\varphi $ and
$\psi $. Then, there are the following possibilities for $G=\Gal
(L/K)$ and $\End _K(A)$:

\begin{enumerate}

\item
$G$ is trivial and $\End _K(A)=\M _0(N)$.
\item
$G=C_2$ and $\End _K^0(A)=\Q (\sqrt {-N})$, $\Q \times \Q $ or $\Q
(\sqrt {m})$ for $m>1$, $m\mid N$, such that $\M _2(\Q )=(\frac
{-N, m}{\Q })$.

\item
$G=C_4$ and $\End _K^0(A)=\Q (\sqrt {-1})$.

\item
$G=D_2$ or $D_4$ and $\End _K^0(A)=\Q $.

\end{enumerate}
\end{theorem}

A third aspect that we will regard concerning the arithmetic of
abelian surfaces with QM stems from the following result obtained
independently by M.\ Jacobson and M.\ Ohta.

\begin{theorem}\cite{Ja}, \cite{Oh}
\label{Oh}
Let $A/K$ be an abelian surface with quaternionic multiplication
by a maximal order $\cO $ in a division quaternion algebra over an
extension $L$ of $K$. Let $ \{ \sigma_\ell \}$ be the compatible
family of Galois representations given by the action of
$\Gal(\qbar / K)$ on the Tate modules $T_\ell(A)$ of $A$ and let
$H=\Gal(\qbar/L)$. Then $\sigma_\ell|_H = \rho_\ell \oplus \rho_\ell$
with
$$
\rho_\ell : H \rightarrow \Aut_\cO (T_\ell(A)) \simeq \cO_\ell^*
$$
and $\rho_\ell$ is surjective for almost every prime.
\end{theorem}

We will obtain an explicit description of the action of the
absolute Galois group $G_K=\Gal(\qbar/K)$ on the Tate modules of
abelian surfaces $A$ with quaternionic multiplication. This allows
us to characterize the three possibilities for $\Gal(L/K)$
described in theorem \ref{L}, $\bf{I}$, and we show how to
effectively determine the field extension $L/K$. Moreover, we
explain how to explicitly bound the finite set of {\em exceptional
primes}, those where the surjectivity conclusion in
Jacobson-Ohta's theorem fails, and illustrate it in a concrete
example with $\Gal(L/K)=C_2$. For all non-exceptional primes, we
also describe the image of the Galois representation $\sigma_\ell$
in $\GL_4(\Z_{\ell })$.

The paper is organized as follows. We devote next two sections to
study the action of $G_K$ on the ring of endomorphisms $\End
_{\qbar }(A)$ and the N\'{e}ron-Severi $\NS (A_{\qbar })$ group of
$A$ respectively. The combination of the description of both
Galois representations eventually yields the proof of our main
theorem \ref{LL} and, as an immediate consequence, proves theorem
\ref{L}.

In section 4 we study the action of Galois on the Tate modules in
the case that $L/K$ is a quadratic extension. Under this
assumption, we show that the Galois representations behave as in
the case of a modular form with (a single) inner twist (cf.
\cite{Ri3}). Following the results of Ribet, we provide sufficient
conditions on a prime $\ell$ for the image of Galois $\rho _{\ell
}(G_L)$ to be {\em as large as possible}.

In section 5, we consider a concrete example of a Jacobian surface
of a curve $C/\Q({\sqrt{-3}})$ of genus $2$ with maximal
quaternionic multiplication borrowed from \cite{HaMu}. Firstly, we
describe the image of the inertia
subgroup at $\ell$ for the residual $\mod \; \ell$ Galois
representations and then, we give a result (lemma
\ref{teo:irreducible}) to distinguish the cases $\Gal(L/K)= C_2$
and $\Gal(L/K) = D_2$ of theorem \ref{L}, $\bf{I}$. By these
means, we show that, in our example, $L/K$ is a quadratic
extension and explicitly determine the field $L$. Then, we proceed
to the determination, by successive elimination of some special
cases, of the images of Galois, following the ideas of \cite{Ri3}.
The main difference with the techniques used in \cite{Ri3}
(described algorithmically in \cite{DiVi}) is that we are dealing
with representations of the Galois group of a number field $K \neq
\Q$.

\section{The action of $\Gal ({\qbar /K})$ on the endomorphism ring.}
\label{end}

In this section we use Chinburg-Friedman's recent classification
of the finite subgroups of maximal arithmetic Kleinian groups
(\cite{ChFr}) to describe the field of definition of the
quaternionic multiplication on an (unpolarized) abelian surface.

Let $A/K$ be an abelian surface over a number field $K$ with
quaternionic multiplication by an order $\cO $ in an indefinite
quaternion algebra $B$ over $\Q$.

The absolute Galois group $G_K=\Gal (\qbar /K)$ of $K$ acts in a
natural way on the full ring of endomorphisms $\End _{\qbar
}(A)=\cO $ of $A$ that we already identify with $\cO $ and induces
a Galois representation

$$
\gamma :G_K\lra \Aut (\cO ).
$$
The Skolem-Noether theorem (\cite{Vi}) asserts that all
automorphisms of a quaternion algebra are inner. Therefore $\Aut
(B)\simeq B^*/\Q ^*$ and the group of automorphisms of $\cO $ is
$\mathrm{N} _{B^*}(\cO )/\Q ^*$ where $\Norm _{B^*}(\cO )=\{
\gamma \in B^*, \gamma ^{-1}\cO \gamma \}$ is the normalizer group
of $\cO $. For $\tau \in G_K$ we will denote $[\gamma _{\tau
}]:B\ra B$ the automorphism of $B$ such that $\beta ^{\tau } =
\gamma _{\tau }^{-1}\beta \gamma _{\tau }$ for any $\beta \in \End
_{\qbar }(A)=\cO $.

If we let $L/K$ be the minimal (and hence normal by proposition
\ref{teo:SiRi}) field extension of $K$ such that $\End _{\qbar
}(A)=\End _L(A)=\cO $, we obtain an exact sequence of groups

$$
1\ra G_L\ra G_K\ra \Norm _{B^*}(\cO )/\Q ^*
$$
and thus a monomorphism $\Gal (L/K)\hookrightarrow \Norm
_{B^*}(\cO )/\Q ^*$.

\begin{proposition}
\label{ChFr}

Let $A/K$ be an abelian surface with quaternionic multiplication
by an order $\cO $, $D=\disc (\cO )$, and let $L/K$ be the minimal
extension of $K$ over which all endomorphisms of $A$ are defined.
Then $L/K$ is either cyclic or dihedral with $\Gal (L/K)\simeq
C_n$ or $D_n$, $n=1, 2$, $3$, $4$ or $6$.

\begin{enumerate}
\item
If $\Gal (L/K)\simeq D_2$, then $B=( \dfrac {d, m}{\Q })$ for some
$d, m\in \Z $, $d, m\mid D$.

\item
If $\Gal (L/K)\simeq C_3$, then any ramified prime $p\mid D$,
$p\not = 2, 3$, satisfies $p\equiv -1 (\mod 3)$; if $\Gal
(L/K)\simeq D_3$, then in addition $B=( \dfrac {-3, m}{\Q })$ with
$m\mid D$.

\item
If $\Gal (L/K)\simeq C_4$, then $2\mid D$ and any odd ramified
prime $p\mid D$ satisfies $p\equiv -1 (\mod 4)$; if $\Gal
(L/K)\simeq D_4$, then in addition $B=( \dfrac {-1, m}{\Q })$ with
$m\mid D$.

\item
If $\Gal (L/K)\simeq C_6$, then $3\mid D$ and any $p\mid D$,
$p\not = 2, 3$, satisfies $p\equiv -1 (\mod 3)$; if $\Gal
(L/K)\simeq D_6$, then in addition $B=( \dfrac {-3, m}{\Q })$ with
$m\mid D$.
\end{enumerate}

\end{proposition}

{\em Proof: } As we already observed, $\Gal (L/K)$ is a finite
subgroup of $\Norm _{B^*}(\cO )/\Q ^*$. Let us firstly note that,
if $\gamma \in B^*$ normalizes the order $\cO $, then $\n (\gamma
)\in \Q ^*$ has odd $p$-adic valuation at any non-ramified prime
number $p\nmid D$. This holds because $\Norm _{\GL _2(\Q _p)}(\M
_2(\Z _p))= \Q _p^* \GL _2(\Z _p)$.

In \cite{ChFr} $\S 2$, Chinburg and Friedman proved that the only
possible finite subgroups of $B^*/\Q ^*$ are the cyclic groups
$C_n$, the dihedral groups $D_n$ and $S_4$, $A_4$ and $A_5$.

By \cite{ChFr}, lemma 2.8, a necessary condition for $B^*/\Q ^*$
to contain either $S_4$, $A_4$ or $A_5$ is that $B=(\frac {-1,
-1}{\Q })$ and this can not be the case because $B$ is indefinite.

Lemma 2.1 in \cite{ChFr} yields that $B^*/\Q ^*$ contains a cyclic
group of order $n>2$ if and only if there exists $\zeta _n\in B^*$
satisfying $\zeta _n ^n = 1$, $\zeta _n ^{n/d}\not =1$ for any
proper divisor $d$ of $n$. In this case, any subgroup
$C_n\subseteq B^*/k^*$ is conjugated to $\langle [1+\zeta
_n]\rangle $. In our case, $\zeta _n\in B^*$ generates a quadratic
field extension $\Q (\zeta _n)/\Q $ and this is only possible for
$n=3, 4$ and $6$. In addition, since $\n (1+\zeta _n)=1, 2$ and
$3$ respectively, the condition $1+\zeta _n\in \Norm _{B^*}(\cO )$
implies that $2\mid D$ for $n=4$ and $3\mid D$ for $n=6$.

It follows from \cite{ChFr}, lemma 2.2, that the conjugacy classes
of subgroups of $\Norm _{B^*}(\cO )/\Q ^*$ of order two are in
bijection with the set of divisors $m\mid D$ of $D$, $m\not =1$,
such that $p$ does not split in $\Q (\sqrt {m})$ for any prime
$p\mid D$, together with $m=1$ if $B\simeq \M _2(\Q )$. This set
is always non-trivial because at least $\pm D$ satisfy these
conditions.

Finally, Chinburg and Friedman proved that $B^*/\Q ^*$ contains a
dihedral subgroup $D_n$, $n\geq 2$, if and only if it contains a
cyclic group $C_n$ (\cite{ChFr}, lemma 2.3). If $n=2$, any
subgroup of $B^*/\Q ^*$ isomorphic to $D_2=C_2\times C_2$ is of
the form $\langle [x], [y]\rangle \subset B^*/\Q ^*$ with $x, y\in
B^*$, $x^2=d$, $y^2=m$, $x y = -y x$ for some $d, m\in \Q ^*$. It
follows that $\Norm _{B^*}(\cO )/\Q ^*$ contains a dihedral group
$D_2$ if and only if $B=( \dfrac {d, m}{\Q })$ for some $d, m\in
\Z $, $d, m\mid D$. Similarly, if $n=3, 4$ or $6$, $\Norm
_{B^*}(\cO )/\Q ^*$ contains a dihedral subgroup $D_n$ if and only
if $B=( \dfrac {d, m}{\Q })$ with $d=-1$ if $n=4$, $d=-3$ if $n=3$
or $6$ and $m\in \Z $, $m\mid D$. In this case $D_n= \langle
[1+\zeta _n], [y]\rangle \subset B^*/\Q ^*$ for some $y\in B^*$,
$y^2=m$. $\Box $

\section{The action of $\Gal (\qbar /K)$ on the N\'{e}ron-Severi group}

Let $A$ be an abelian variety defined over a number field $K$. For
any field extension $L/K$, we let $A_L = A \times _K L$ denote the
same abelian variety $A$ with the base extended to Spec $L$. Let
$\Div (A)$ denote the group of Weil divisors of $A$ and
let $\Pic (A)$ denote the group of invertible sheaves on $A$ over $K$.

Let $\Pic ^0(A_{\qbar })$ denote the subgroup of $\Pic (A_{\qbar
})$ of invertible sheaves algebraically equivalent to $0$ and let $\Pic
^0(A) = \Pic (A)\cap \Pic ^0(A_{\qbar })$. The N\'{e}ron-Severi group $\NS (A)$ of $A$ is $\NS (A) = \Pic
(A)/\Pic ^0(A)$. The algebraic class of an invertible sheave $\cL $ lies in $H^0(\NS (A_{\qbar }))=H^0(\Gal (\qbar /K), \NS (A_{\qbar }))$ if and only if all its Galois
conjugates $\cL ^{\tau }$, $\tau \in \Gal (\qbar /K)$, are algebraically
equivalent to $\cL $. We define the Picard number of $A_K$ to be
$\rho (A_K)= \rank _{\Z } H^0(\NS (A_{\qbar }))$; it is a finite
number due to N\'{e}ron's basis theorem.

Let now $\cO $ be an hereditary order in a quaternion algebra $B$
over $\Q $. Assume that $A$ is an abelian surface defined over a
number field $K$ together with an isomorphism of rings $\iota :\cO
\stackrel {\sim } {\ra }\End _{\qbar }(A)$. The underlying complex
torus $A_{\C}=V/\Lambda $ is the quotient of a complex vector
space $V$ of dimension $2$ by a lattice $\Lambda $ of rank $4$
over $\Z $. Upon fixing an isomorphism $B\otimes \R \stackrel
{\sim }{\ra } \M _2(\R )$ there is an action of $\cO \subset
B\subset \M _2(\R )$ on the lattice $\Lambda $ that makes it a
left $\cO $-module. Since all left ideals of $\cO $ are principal
and from the work of Shimura (\cite{Sh2}), it is well-known that
there exists $\tau \in \cH =\{ a+b i\in \C , b>0\} $ such that
$\Lambda =\cO \begin{pmatrix}\tau \\ 1
\end{pmatrix}$.

In \cite{Ro1}, the absolute N\'{e}ron-Severi group $\NS (A_{\C })\simeq \NS
(A_{\Q })$ was largely studied under the assumption on $\cO $ to
be hereditary: it was seen that the first Chern class allows us to
regard $\NS (A_{\C })$ as a sub-lattice of the $3$-dimensional
vector space $B_0=\{\mu \in B, \tr (\mu )=0\}$ of pure quaternions
of $B$ in a way that fundamental properties of line bundles $\cL $
on $A$ such as the degree $\deg (\cL )$ (\cite{Ro1}, proposition
3.1), the behaviour under pull-backs by endomorphisms (\cite{Ro1},
theorem 2.2) and the index $i(\cL )$ and the ampleness
(\cite{Ro1}, theorem 5.1) can be interpreted in terms of the
arithmetic of $B$. We summarize it in the following

\begin{theorem}\cite{Ro1}
\label{Ro1} Let $A/\qbar $ be an abelian surface with $\End
_{\qbar }(A) \stackrel {\iota }{\simeq } \cO $ an hereditary order
of discriminant $D$ in a quaternion algebra. Then there is an
isomorphism of additive groups

$$
\begin{matrix}
c_1: & \NS (A_{\qbar }) & \ra & \cO _0^{\sharp } \\
     & \cL             & \mapsto & c_1(\cL ) \\
\end{matrix}
$$
such that
\begin{enumerate}
\item
$\deg (\cL ) = D\cdot \n (c_1(\cL ))$.

\item
For any endomorphism $\alpha \in \cO =\End _{\qbar }(A)$,
$c_1(\alpha ^*(\cL )) = \bar {\alpha }c_1(\cL )\alpha $.

\item
A line bundle $\cL \in \NS (A_{\qbar })$ is a polarization if and
only if $\n (c_1(\cL ))>0$ and $\det (\nu _{\cL })>0$ where $\nu
_{\cL }\in \GL _2(\R )$ is (any) matrix such that $\nu _{\cL
}^{-1} c_1(\cL ) \nu _{\cL }\in \Q ^* \begin{pmatrix} 0 & 1 \\ -1
& 0 \end{pmatrix}$.

\end{enumerate}
\end{theorem}
Here, $\cO ^{\sharp }=\{ \beta \in B, \tr (\cO \beta )\subseteq \Z
\}$ denotes the codifferent ideal of $\cO $ in $B$. By $\cO _0
^{\sharp }$ we mean the subgroup $\cO ^{\sharp }\cap B_0$ of pure
quaternions of $\cO ^{\sharp }$. For our purposes in this note, we
only need to know that it is a lattice in $B_0$ and in particular
$\rho (A_{\qbar })=3$.

Let us also remark that, by Eichler's theory on optimal embeddings
(cf. e. g. \cite{AlBa}, \cite{Vi}), there always exists $\mu \in
\cO $ such that $\mu ^2+D=0$ and, as a corollary of theorem
\ref{Ro1}, $A_{\qbar }$ is always principally polarizable. We
refer the reader to \cite{Ro1} for more details.

We consider now the action of the Galois group $G_K = \Gal (\qbar
/K)$ on $\NS (A_{\qbar })$ given by $\cL (\cD )^{\tau } = \cL (\cD
^{\tau })$ for any line bundle $\cL $ on $A$ represented by a Weil
divisor $\cD $ and $\tau \in G_K$. From theorem \ref{Ro1}, any
automorphism of $\NS (A_{\qbar })$ can be regarded as a linear
automorphism of $B_0$. Moreover, since the Galois action preserves
the degree of line bundles and the first Chern class is a
monomorphism of quadratic modules $c_1:(\NS (A_\qbar ), \deg
)\hookrightarrow (B_0, D\cdot \n )$, we obtain a Galois
representation

$$
\begin{matrix}
\eta : & G_K & \lra & \Aut (\NS (A_{\qbar }), \deg )&\subset &\Aut (B_0, D\cdot \n )\\
       & \tau &\mapsto & & &\eta _{\tau }\\
\end{matrix}
$$

We have that

\begin{enumerate}
\item
For any $\alpha \in \cO =\End _{\qbar }(A)$, $(\alpha ^*(\cL )^{\tau }) =
(\alpha ^{\tau })^*(\cL ^{\tau })$.

\item
The index $i(\cL )$ only depends on the $G_K$-orbit of $\cL $,
that is, $i(\cL ^{\tau})=i(\cL )$, for any $\tau \in G_K$. In
particular $\cL ^{\tau }$ is a polarization if and only if $\cL $
is.

\end{enumerate}

The following relates the Galois actions on $\End _{\qbar }(A)$
and on the N\'{e}ron-Severi group of $A/K$ by means of a reciprocity
law.

\begin{theorem}
\label{law}

Let $A/K$ be an abelian surface with QM by an hereditary order
$\cO $ of discriminant $D$ in a quaternion algebra $B$.

Let $\gamma :G_K\ra \Aut (\End _{\qbar }(A))$, $\tau \mapsto
[\gtau ]: \cO \ra \cO , \beta \mapsto \gtau ^{-1}\beta \gtau $, be
the action of $\Gal (\qbar /K)$ on the ring of endomorphisms of
$A$. Define $\etau  = \sign (\n (\gtau ))=\pm 1$. Then

$$
c_1(\cL ^{\tau }) = \etau \cdot \gtau ^{-1} c_1(\cL ) \gtau
$$
for any line bundle $\cL \in \NS _{\qbar }(A)$ and any $\tau \in
G_K$.

\end{theorem}

{\em Proof: } Fix $\tau \in G_K$. We firstly claim that $\eta
_{\tau }:$ $B_0\ra $ $B_0$ is given by $\mu \mapsto \et \gt
^{-1}\mu \gt $ for some $\et =\pm 1$, $\gt \in B^*$. Indeed, any
linear endomorphism of $B_0$ extends uniquely to an endomorphism
of $B$ and $\End (B)\simeq B\otimes B$ with $\gt _1\otimes \gt
_2:B\ra B$, $\beta \mapsto \gt _1\beta \gt _2$. We must have in
addition that $\tr (\gt _1 \mu \gt _2)=\tr (\gt _2 \gt _1 \mu )=0$
for any $\mu \in B_0$. This automatically implies that $\gt _2 \gt
_1\in \Q $.

Since the action of $G_K$ on $\NS _{\qbar }(A)$ conserves the
degree of line bundles, we deduce from theorem \ref{Ro1} that $\n
(\eta _{\tau }(\mu )) = \n (\gt _1 \mu \gt _2) =\n (\gt _1)\n (\mu
)\n (\gt _2)= \n (\mu )$ for any $\mu \in B_0$. Hence $\n (\gt
_2)= \n (\gt _1)^{-1}$ and thus $\gt :=\gt _2 = \et \gt _1^{-1}$
for some $\et =\pm 1$. This proves the claim.

We now show that $\gt = \gtau \in B^*/\Q ^*$ and $\et =\etau $. We
know that $(\alpha ^*(\cL )^{\tau }) = (\alpha ^{\tau })^*(\cL
^{\tau })$ for any $\alpha \in \cO $. Taking theorem \ref{Ro1}
into account this implies that $\eta _{\tau }(\bar {\alpha }\mu
\alpha ) = \bar {[\gtau ](\alpha )} \eta _{\tau }(\mu ) [\gtau
](\alpha )$ and thus $\et \gt ^{-1}(\bar {\alpha } \mu \alpha )
\gt  = \et (\bar {\gtau ^{-1}\alpha \gtau }) \gt ^{-1}\mu \gt
(\gtau ^{-1}\alpha \gtau ) $ for any $\alpha \in B$, $\mu \in $
$B_0$. Choosing $\alpha = \mu $ and bearing in mind that $\gtau
^{-1}=\bar {\gtau }\n (\gtau )^{-1}$, this says that $\gt ^{-1}
\mu \gt = \gtau ^{-1} \mu ^{-1} \gtau \gt ^{-1} \mu \gt \gtau
^{-1} \mu \gtau $ and thus $\mu (\om \mu \om ^{-1})=(\om \mu \om
^{-1}) \mu $, where we write $\om = \gtau \gt ^{-1}$. The
centralizer of $\Q (\mu )$ in $B$ is $\Q (\mu )$ itself and
therefore $(\om \mu \om ^{-1})\in \Q (\mu )$. But $\tr (\mu )=\tr
(\om \mu \om ^{-1})=0, \n(\mu )=\n (\om \mu \om ^{-1})$ and this
implies that $\mu = \pm \om \mu \om ^{-1}$. Since this must hold
for any $\mu \in B_0$, it follows that $\om \in \Q ^*$ and thus
$\gt = \gtau \in B^*/\Q ^*$ as we wished.

We then already have that $\eta _{\tau }:$ $B_0\ra B_0$ is given
by $\mu \mapsto \et \gtau ^{-1}\mu \gtau $ for some $\et =\pm 1$.
If $\mu = c_1(\cL )$ for a polarization $\cL $ on $A$, this means
that $c_1(\cL ^{\tau })= \et \gtau ^{-1}\mu \gtau $. Since $\cL
^{\tau }$ is still an ample line bundle we have, according to
theorem \ref{Ro1}, 3., that $\et = \sign (\n (\gtau ))$. $\Box $

We can now prove the following result that yields theorem \ref{L}
in the introduction as an immediate corollary. Let us before say a
word about the N\'{e}ron-Severi group and the cone of polarizations on
an abelian variety $A/K$.

Firstly, we remark that $A_K$ is always attached with a
polarization and thus the Picard number $\rho (A_K)$ never
vanishes. Namely, if $A\stackrel{i}{\hookrightarrow }\PP ^N_K$ is
an embedding of $A$ into a projective space over $K$, then
$i^*(\cO (1))\in \NS (A_K)$. However, any polarization constructed
by these means is very ample and hence can not be principal. This
should be taken into account together with the fact that any
abelian surface with hereditary quaternionic multiplication admits
a principal polarization over $\qbar $ (cf.\,\cite{Ro1}, Corollary
6.3).

By \cite{Mu}, or also theorem \ref{Ro1}, we thus have that $1\leq
\rho (A_K)\leq 3$. Both three cases are possible and each
possibility has a direct translation in terms of the algebra of
endomorphisms: $\rho (A_K)=1$ if and only if $\End _K^0(A):=\End
_K(A)\otimes \Q $ is $\Q $ or an imaginary quadratic field, $\rho
(A_K)=2$ if and only if $\End _K^0(A)$ is a real quadratic field
and $\rho (A_K)=3$ if and only if $\End _K^0(A)=\End _{\qbar }(A)=B$.

The following definition was introduced in \cite{Ro2} and
\cite{Ro3}.
\begin{definition}

Let $\cO$ be an order in a quaternion algebra $B$ over $\Q $ and
let $D=\disc (\cO )$. We say that $\cO $ admits a {\em twist of
degree} $\delta \geq 1$ if there exists $m\in \Z $, $m\mid D$,
such that

$$
B = \Q +\Q i + \Q j + \Q i j = (\dfrac {-D \delta , m}{\Q })
$$
with $i$, $j\in \cO $, $i^2=-D \delta $, $j^2=m$ and $i j = -j i$.
In this case we say that the twist is of {\em norm} $m$. If
$\delta =1$, we say that $\cO $ admits a {\em principal twist}.

\end{definition}

{\em Remark.} For a given order
$\cO $ in an indefinite quaternion algebra, let
$\mathcal {N}_{\delta }=\{ m_1, ..., m_t\} $,
$0<m_i\mid D$, denote the (possibly empty) set of norms of the
twists of degree $\delta
$ on $\cO $. It is easy to show that $\mathcal {N}_1$ is either empty or
$\mathcal {N}_1=\{ m, D/m\} $ for some
$m\mid D$. In other cases, $\mathcal {N}_{\delta }$ can be larger.
Indeed, if $\delta =D$ for instance, then $\mathcal {N}_D$ is either
empty or equal to the set of sums of two squares $m=m_1^2+m_2^2$ that divide $D$. We
finally note that a
quaternion order $\cO $ can very well admit twists of several different
degrees.

{\em Remark. } In practice, the
computation of a finite number of Hilbert symbols suffices to decide whether a
given indefinite order is
twisting of certain degree $\delta $. Let us just
quote that a necessary and sufficient condition for $B$ to contain
a maximal order $\cO $ admitting a twist of degree $\delta
$ and norm $m$ is that $m>0$, $m\mid D=\disc (\cO )=\disc (B)$ and
that for any odd prime $p\mid D$: $m\not \in \F _p^{*^2}$ if $p\nmid m$
($D/m\not \in \F _p^{*^2}$ if $p\mid m$ respectively).

Examples of quaternion algebras with principally twisting maximal
orders are $B = (\dfrac {-6, 2}{\Q }  )$ and $B = (\dfrac {-10,
2}{\Q } )$ of discriminant $D = 6$ and $D = 10$ respectively.

\begin{theorem}
\label{LL}

Let $A/K$ be an abelian surface defined over a number field $K$
with quaternionic multiplication by an hereditary order $\cO $ of
discriminant $D$ in a quaternion algebra $B$ and let $L/K$ be the
minimal extension of $K$ such that $\End _L(A)\simeq \cO $. Fix a
polarization $\cL _0\in H^0(\Gal (\qbar ,K),\NS (A_{\qbar }))$ and let $\delta =\deg (\cL _0)$ be
its degree.

{\bfseries A.}
\begin{enumerate}
\item
If $\delta $ is not equal to $D$ neither to $3 D$ up to squares,
then $\Gal (L/K) \simeq \{ 1 \}$, $C_2$ or $D_2 = C_2\times C_2$.

\item
If $\delta = D k^2$ for some $k\in \Z $, then $\Gal (L/K)\simeq
C_n$ or $D_n$ with $n = 1$, $2$ or $4$.

\item
If $\delta = \frac {D k^2}{3}$ for some $k\in \Z $, then $\Gal
(L/K)\simeq C_n$ or $D_n$ for $n = 1$, $2$, $3$ or $6$.

\end{enumerate}

{\bfseries B.} In any of the cases above, if $\cO $ does not admit
any twist of degree $\delta $, then $\Gal (L/K)$ is necessarily
cyclic.

{\bfseries C.}
\begin{enumerate}
\item
If $\Gal (L/K)\simeq C_2$, then $\End _K^0(A)\simeq \Q (\sqrt {-D
\delta })$ or $\Q (\sqrt {m_i})$ for $0<m_i\in \mathcal
{N}_{\delta }$ a norm of a twist of degree $\delta $ on $\cO $.

\item
If $G=C_3$ or $C_6$, then $\End _K^0(A)\simeq \Q (\sqrt {-3})$.

\item
If $G=C_4$, then $\End _K^0(A)\simeq \Q (\sqrt {-1})$

\item
If $G=D_n$, then $\End _K(A) \simeq \Q $.

\end{enumerate}

\end{theorem}

{\em Proof: } Recall that, according to proposition \ref{ChFr},
$\Gal (L/K) \simeq C_n$ or $D_n$ with $n=1, 2$, $3$, $4$ or $6$. Let
$\cL _0\in H^0(\Gal (\qbar ,K),\NS (A_{\qbar }))$ be a polarization on $A$ and let $\mu = c_1(\cL _0)\in
\cO _0$. It satisfies that $\mu ^2 + D \delta  = 0$ by theorem
\ref{Ro1}, (1). Fix $\tau \in \Gal (L/K)$ and let
$\gtau \in B^*$ the quaternion associated to $\tau $ in section 2.
Suitably scaling it, we can (and we do) choose a representative in
$\Norm _{B^*}(\cO )/\Q ^*$ such that $\gtau \in \cO $ and $\n
(\gtau )$ is a square-free integer. Then, since $\gtau $ must
normalize $\cO $, we know that $\n (\gtau )\mid D$.

Since the algebraic class of $\cL _0$ is $\Gal (\qbar /K)$-invariant, it
follows from theorem \ref{law} that $\mu =c_1(\cL _0) = c_1(\cL
_0^{\tau }) = \etau \gtau  ^{-1} \mu \gtau $.

If $\etau =-1$, then the above expression implies that $\mu \gtau
= -\gtau \mu $. Since $\tr (\mu \gtau )=\mu \gtau - \bar {\gtau
}\mu =-\tr (\gtau ) \mu \in \Q ^*$, we deduce that $\tr (\gtau
)=0$. This means that $\gtau ^2=m$ for some $m\mid D$ and $B = \Q
+ \Q \mu + \Q \gtau + \Q \mu \gtau = (\frac {-D \delta , m}{\Q
})$. The indefiniteness of $B$ forces $m$ to be positive. We
obtain that in this case $\langle [\gtau ] \rangle \simeq C_2$.

On the other hand, if $\etau =1$, then $\gtau \in \Q (\mu )\simeq
\Q (\sqrt {-D \delta })=\Q (\sqrt {-\underline {D \delta }})$,
where we let $\underline {D \delta }$ denote the square-free part
of $D \delta $. In this case, and bearing in mind that $\gtau $
must generate a finite subgroup of $B^*/\Q ^*$, we deduce that
either

\begin{itemize}
\item
$\gtau = \pm \sqrt {\underline {D \delta }/D \delta } \cdot \mu $
and hence $\underline {D \delta }\mid D $ and $\langle [\gtau
]\rangle \simeq C_2$,

\item
$\gtau = 1+\zeta _n$ for some $n^{\mbox{th}}$-primitive root of
unity $\zeta _n\in B^*$, $n=3$ or $6$, and hence $\underline {D
\delta } = 3$ and $\langle [\gtau ]\rangle \simeq C_n$ or

\item
$\gtau = 1+\zeta _4$ for some $4^{\mbox{th}}$-primitive root of
unity $\zeta _4\in B^*$ and hence $\underline {D \delta } = 1$ and
$\langle [\gtau ]\rangle \simeq C_4$.

\end{itemize}
We conclude that a necessary condition for $\Gal (L/K)$ to contain
a cyclic subgroup of order $n\geq 3$ is $\underline{D \delta } =
1$ or $3$ which amounts to say that $\deg (\cL _0)=\delta $ is $D$
or $3 D$ up to squares respectively. Also, if $\deg (\cL _0)= D
k^2$, then necessarily $\Gal (L/K)\simeq C_n$ or $D_n$ with $n=1$,
$2$ or $4$ and an analogous statement holds if $\deg (\cL _0)= 3
D$ up to squares. Further, if $B\not \simeq (\frac {-D \delta ,
m}{\Q })$ for any $0<m\mid D$, then it follows from the discussion
above that $\epsilon _{\tau }=1$ for any $\tau \in \Gal (L/K)$
and, as a consequence, $\Gal (L/K)\subset \Q (\mu )^*/\Q ^*$.
Since the only finite subgroups of $\Q (\mu )^*/\Q ^*$ are cyclic,
the proof of parts {\bfseries A} and {\bfseries B} is completed.

As for part {\bfseries C}, assume first that $\Gal (L/K)=\langle
[\gtau ] \rangle \simeq C_2$. Then $\gtau \in B^*$ satisfies
$\gtau ^2=-\n (\gtau )\in \Q ^*$ and we already saw that the only
possibilities are, up to squares, $\n (\gtau ) = D \delta $ or
$m\in \mathcal {N}_ {\delta }$. In any of these cases, $\End
_K^0(A)=\{ \beta \in \End _L(A): \beta ^{\tau }=\beta \} = \{
\beta \in \End _L(A): \beta \gtau =\gtau \beta \} = \Q (\gtau )$
and this implies our first assertion of part {\bfseries C}.
Similarly, if $\Gal (L/K)=\langle 1+\zeta _n\rangle \simeq C_n$
with $n=3$, $4$ or $6$ then $\End _K^0(A)=\Q (1+\zeta _n)\simeq \Q
(\sqrt {-1})$ or $\Q (\sqrt {-3})$ depending on the cases.
Finally, if $\Gal (L/K)=\langle \gtau , \gamma _{\tau '}\rangle
\simeq D_n$ with $\langle \gtau \rangle \simeq C_n$ and $\langle
\gamma _{\tau '}\rangle \simeq C_2$, then $\End _K^0(A)=\{ \beta
\in \Q (\gtau ): \beta ^{\tau '}=\beta \} = \Q $. Here, the last
equality holds because it is not possible that $\gtau $ and
$\gamma _{\tau '}$ commute. $\Box $

The following lemma may be useful in many situations in order to
apply theorem \ref{LL}. It easily follows from proposition
\ref{Ro1}.

\begin{lemma}

Let $A/\qbar $ be an abelian surface with $\End (A)$ a maximal
order in a quaternion algebra of discriminant $D$. If there exist
prime numbers $p, q\mid D$ such that $p$ splits in $\Q (\sqrt
{-1})$ and $q$ splits in $\Q (\sqrt {-3})$, then no polarizations
on $A$ have degree $D k^2$ or $3 D k^2$ for any $k\in \Z $.

\end{lemma}

\section{The action of $\Gal(\qbar /K)$ on Tate modules}

Let $A$ be an abelian surface defined over a number field $K$ such
that $\End _{\qbar }(A) \simeq \cO $ is an order in an indefinite
quaternion algebra $B$ over $\Q$. Let $L/K$ be the minimal field
of definition of the endomorphisms of $A$. In this section we
consider the compatible family of Galois representations $\{
\sigma_\ell \}$ given by the action of $G= \Gal ( \qbar /K )$ on
the Tate modules $T_\ell(A)$ of A. Throughout, we restrict
ourselves to primes $\ell \nmid D \cdot N$, where $D = \disc(\cO
)$ and $N$ is the product of the places of bad reduction of $A$.
We refer the reader to \cite{Jor} for an accurate study of the
Galois representations arising from ramified primes $\ell \mid
\disc (B)\mid D$.

Our aim is two-fold: we wish to make an effective approach to
Jacobson-Ohta's theorem \ref{Oh} and to our main result \ref{LL}.
Bearing this idea in mind, we assume that $L/K$ is a quadratic
extension, in contraposition to the other possibilities permitted
by theorem \ref{LL}. Then, we know that $\End_K(A)$ is an order in
a quadratic field $\Q(\sqrt{d})$. Assume further for simplicity
that $\Q (\sqrt {d})$ is real, and let $\cR$ be its ring of
integers.
 The imaginary quadratic case is
described along the same lines but for certain differences (cf.
section $5.2$).

In this setting, it is well-known that the four dimensional Galois
representations $\sigma _{\ell }$ are reducible over
$\Q(\sqrt{d})$, that is,

$$
\sigma_\ell = \re \oplus \re^\gamma ,
$$
where $\la \mid \ell$ is a prime in $\Q(\sqrt{d})$ over $\ell $
and $\Gal(\Q(\sqrt{d})/\Q)=\langle \gamma \rangle $. Moreover,
$\det(\re) = \chi $ is the $\ell$-adic cyclotomic character (cf.
\cite{Ri4}).

The representations $\re$ have their images contained in the
groups

$$
\U = \{ M \in \GL_2(\cR_\la) \; \mid  \; \det M \in \Z_\ell^* \} .
$$

If we consider the subgroup $H= \Gal(\qbar/L)$ of index two of
$G$, we know by theorem \ref{Oh} that $\re|_H : H \rightarrow
\Aut_\cO(T_\ell(A)) \simeq \cO_\ell^*\simeq \GL_2(\Z_\ell)$ is
surjective for almost every prime $\ell$. This imposes a strong
restriction on the image of $\re$: in general it can not be the
full group $\U$ but a subgroup $\V := \Image(\re )$ that contains
the image of $H$ as a normal subgroup of index at most $2$.

For a prime $\la $, let us call $H_\ell := \Image (\re|_H)$. We
know that it always holds that $H_\ell \subseteq \Ol$ and we say
that $\ell $ is an {\em exceptional} prime if the inclusion fails
to be an equality.

By applying Faltings' proof of Tate's Conjecture as in \cite{Ri4}, Prop. 3.5 (together with Cebotarev density theorem)
the condition on $\End_K(A) \subseteq \cR $ to be a real quadratic order
implies that there are infinitely many Frobenius elements with
$\tr (\re(\Fr \; \wp)) \in \cR \smallsetminus \Z$. In turn, this
implies that $\V \nsubseteq \Ol $ and that $H_{\ell }$ is a normal
subgroup of $\V$ of index two for almost every prime $\lambda$
such that $(  \frac{d}{\ell} ) = -1$. If, on the other hand,
$(\frac{d}{\ell}) = 1$, then $\cR_\la = \Z_\ell$ and thus $\V
\subseteq \Ol$. In particular, by applying theorem \ref{Oh}, it
holds that $\Ol \subseteq \V \subseteq \U $ and $[\V : \Ol ]= 2$
if $\left (\dfrac {d}{\ell }\right ) = -1$ and $\Ol = \V $ if
$\left (\dfrac {d}{\ell }\right ) = 1$, for almost every prime
$\la $.

Let us restrict for a while to inert primes: those such that
$(\frac{d}{\ell}) = -1$. For them, we have that $\V \subseteq \U$
and there is the following exact sequence

$$
0 \rightarrow H_\ell \rightarrow \V \rightarrow \{ \pm 1 \}
\rightarrow 0 \quad \quad \quad \quad \eqno(4.1)
$$

Moreover, it is easy to see that $H_\ell = \V \cap \Ol$. By
considering the quotient $\V / H_\ell \simeq \Gal(L/K)$, using the
information on the ramification of $\re$ and varying the prime
$\ell$, we conclude that $L/K$ is unramified outside $N$ in
agreement with proposition \ref{teo:SiRi}. From (4.1), it follows
that the quadratic character $\psi $ corresponding to $\Gal(L/K)$
determines whether or not the image $\re(\Fr \; \wp)$ belongs to
$H_\ell$, i.e., $\re (\Fr \; \wp ) \in H_\ell \Longleftrightarrow
\psi(\wp)= 1$. Equivalently,

$$
\re (\Fr \; \wp ) \in \Ol \Longleftrightarrow \psi(\wp)= 1 \quad
\quad \quad \quad \eqno(4.2)
$$

Let $\wp$ be a prime such that $\psi(\wp) = -1$, so that $\re (\Fr
\; \wp) \in \V \smallsetminus H_\ell$. Since $[\V : H_\ell] =2$,
we have that

$$
\re^2 (\Fr \; \wp)\in H_\ell \subseteq \Ol .\quad \quad \quad
\quad \eqno(4.3)
$$

Let us denote $a_\wp := \tr (\re (\Fr \; \wp))$. We know that the
determinant of $\re$ is the $\ell$-adic cyclotomic character
$\chi$ and we hence obtain that $\tr (\re^2 (\Fr \; \wp) ) =
a_\wp^2 - 2 p $. Thus, (4.3) implies that

$$
a_\wp^2 \in \Z_\ell \quad \quad \quad \quad \eqno(4.4).
$$

Observe that (4.4) is automatic for the remaining primes: those
that split in $\Q(\sqrt{d})/\Q$. Therefore, we have
$\Q(\{a_\wp \})
 = \Q(\sqrt{d})$ and $\Q( \{ a_\wp^2 \}) =\Q $. From this and the fact
that the character $\psi$ governs the behaviour of $\re$ for every
$\ell$ inert in $\Q(\sqrt{d})$ (see (4.2)), it is an easy exercise
to show that, for every $\wp \nmid D\cdot N$:

$$
a_\wp^\gamma = \psi(\wp) a_\wp . \quad \quad \quad \quad \;
\eqno(4.5)
$$

In fact, if the determinant is defined over $\Q $, Serre proved
that compatible families of Galois representations verifying this
property of having {\em inner twists} are characterized  by the
strict inclusion of $\Q(\{ a_\wp^2 \})$ in $\Q( \{ a_\wp \})$.

For an arbitrary rational prime $\ell$, observe that $a_\wp =
u_\wp \sqrt{d}$, $u_\wp \in \Z $ if $\psi(\wp) = -1 $. Let us fix
such an element with the further restriction $a_\wp \neq 0$ and
let $p$ be the rational prime such that $\wp \mid p$. To ease  the
notation, we denote $ M_\wp =
\begin{pmatrix}
  a_\wp & 0 \\
  0 & 1/a_\wp
\end{pmatrix}$.

Using (4.5) and imitating the proof of the theorem of Papier (see
\cite{Ri3}, section 4 or \cite{Di} pg. 398) with the restriction that $\ell \nmid p
\Norm(a_\wp)$, we deduce that $\V = \langle  \GL_2(\Z_\ell) ,
M_\wp \rangle $ for those non-exceptional primes $\ell $. In
particular, this holds for almost every prime $\ell$. In
conclusion, we have shown that for every prime $\ell$ verifying
$$
\ell \nmid D\cdot N \cdot p \cdot \Norm(a_\wp), \quad \quad \quad
\; \eqno(4.6)
$$
the maximal possible image of $\re$ is $\Olt:=  \langle
\GL_2(\Z_\ell) , M_\wp \rangle \subseteq \GL_2(\cR_\lambda)$ and
that the image is in fact maximal for almost every $\lambda$. We
thus wonder: {\em How can the finite set of rational primes $\ell
$ such that $\V \varsubsetneq \Olt $ be bounded?}

In addition to restriction (4.6), the main point is that we need
to impose the following condition on the residual $\mod \;
\lambda$ representations $\rer$ (obtained from $\re$ by composing
with the naive reduction): $\rer|_H$ must be irreducible and the
order of the image $\Image (\rer |_H)$ must be a multiple of
$\ell$.

As explained in \cite{Ri3}, for $\ell > 2 \geq [\V : H_\ell ]$, if
one checks that $\rer$ is irreducible and that the order of its
image is a multiple of $\ell$, then $\rer|_H$ will also verify
both these conditions. The main result on the determination of
images for $2$-dimensional Galois representations in \cite{Ri3}
(see also \cite{Ri1} and \cite{Di}) then implies that $H_\ell =
\GL_2 (\Z_\ell)$ and $\V = \Olt$ for every prime $\la $, $\la \mid
\ell $ such that

\begin{itemize}

\item
$\ell$ verifies (4.6)

\item
$\rer$ is irreducible with image of order multiple of $\ell$

\item
$\ell \geq 5$ and $\det: H_\ell \rightarrow \Z_\ell^*$ is
surjective

\item there exists a prime $\Re$ in $K$ with $a_\Re \in \mathbb{Z}$
and $\ell \nmid a_\Re, \; \Re \nmid \ell$

\end{itemize}

Since $\det(\re|_H)= \chi|_H$, $\det |_{H_{\ell }}$ is surjective
whenever $\ell$ does not ramify in $L/\Q$. As $L/K$ is unramified
outside $N$, it is enough to impose that $\ell \nmid N, \; \ell
\nmid \disc(K)$.

\begin{theorem}
\label{teo:ribetimages}

Let $\wp$, $\Re$ be primes in $K$ such that $\Q(a_\wp) =
\Q(\sqrt{d})$ and $a_\Re \in \Z$. Let $\ell \geq 5$ be a rational
prime such that $\ell \neq p , \; \ell \nmid \Norm(a_\wp), \; \ell
\nmid N \cdot D , \; \ell \nmid \disc(K), \;  \ell \nmid a_\Re, \;
\Re \nmid \ell $ and let $\la \mid \ell$ in $\Q(\sqrt{d})$ be such that $\rer
$ is irreducible and $\ell\mid |\Image (\rer )|$. Then,

$$
\Image(\re) = \langle  \GL_2(\Z_\ell) , M_\wp  \rangle .
$$
\end{theorem}

The condition $\ell \mid |\Image(\rer)|$ can be dealt with by
elimination. Using the classification of maximal subgroups of
$\PGL_2$ over a finite field (of characteristic $\ell$) due to
L.E. Dickson, we know that any irreducible subgroup either has
order multiple of $\ell$ or its projective image falls in one of
the following cases: cyclic, dihedral or {\em small exceptional}
(isomorphic to
 $A_4 , S_4$ or $A_5$). Thus, the above theorem  asserts that,
in order to explicitly bound the set of exceptional primes in a
concrete example, it only remains to bound the set of primes such
that $\rer$ (modulo its centre) is either reducible, cyclic,
dihedral or small exceptional. This will be accomplished in the
following section.

Another key ingredient in the determination of the image of $\re$
is the description of the restriction of $\rer$ to the inertia
subgroup at $\ell$: the determinant of   $\re$ being the
cyclotomic character, we know a priori that both $\re$ and its
residual counterpart $\rer$ necessarily ramify at $\ell$. Thanks
to the results of  Raynaud (cf. \cite{Ra}), we know that one of
the following must hold:

$$
\rer|_{I_\ell} \simeq
\begin{pmatrix}
  1 & * \\
  0 & \chi
\end{pmatrix} \; \; \mbox{or} \; \;
\begin{pmatrix}
  \psi_2 & 0 \\
  0 & \psi_2^\ell
\end{pmatrix},
$$
where $\chi$ denotes the $\mod \; \ell$ fundamental character and
$\psi_2$ a fundamental character of level $2$.

\section{A concrete example of $\mathrm{GL} _2$-type}

\subsection{Fields of definition and endomorphism algebras}

In this section we illustrate our results with an explicit
example. We refer to \cite{DiRo} to further examples of Jacobians
with quaternionic multiplication with different behaviours.

Let $C$ be the smooth projective model of the genus $2$
hyperelliptic curve $Y^2 = \frac{1}{48}X(9075X^4+3025(3+2
\sqrt{-3}) X^3 -6875 X^2 +220 (-3+ 2 \sqrt{-3}) X + 48)$. By
\cite{HaMu}, the ring of endomorphisms of the Jacobian variety of
$C$ over $\qbar $ is a maximal order in the quaternion algebra of
discriminant $D= 10$ over $\Q $. In this section we prove the
following.

\begin{theorem}
\label{teo:C2}

Let $J(C)/K$ be the Jacobian variety of $C$ over $K=\Q (\sqrt
{-3})$. Then, $L= \Q( \sqrt{-3}, \sqrt{-11})$ is the minimal field
of definition of the quaternionic endomorphisms of $J(C)$ and

$$
\End_K(J(C))\otimes_\Z \Q = \Q(\sqrt{5}).
$$
\end{theorem}

Moreover, as it is shown in \cite{HaMu}, there is an isomorphism
of curves $C\stackrel {\simeq }{\ra }C^{\tau }$ between $C$ and
the Galois conjugated curve $C^{\tau }$, where $\tau $ denotes the
non trivial involution of $K$ over $\Q $. Since the isomorphism
lifts to an isomorphism $J(C)\simeq J(C)^{\tau }$ of abelian
varieties, the generalized Shimura-Taniyama-Weil Conjecture
predicts that $J(C)$ should be modular (cf. \cite{Ri4}).

According to theorem \ref{L}, in order to prove that $L/K$ is a
quadratic extension, it suffices to exclude the cases $L=K$ and
$\Gal(L/K) =D_2$.

From the model we have of $C$ we see that its set of primes of bad reduction
 is contained in $\{
2, 3, 5, 7, 11\}$. Thus, we take  $N = 2 \cdot \tres \cdot 5 \cdot 7 \cdot
11$, where $\tres = \sqrt{-3}$ ramifies in $K/\Q $. Let us
consider the Galois representations $\sigma_\ell$ acting on the
Tate modules of $A=J(C)$. In \cite{HaMu}, the  characteristic
polynomials $\Pol _\wp (x)$ of the matrices $\sigma_\ell(\Fr \;
\wp)$ for the first primes $\wp \nmid N$ of $K$ of residue class
degree $1$ were computed and factorized as follows:

$$
\Pol _\wp(x) = (x^2-a_\wp x + p) (x^2-b_\wp x + p) \quad \quad
\quad \quad \eqno(5.1)
$$

The computed values of $a_\wp$, $b_\wp$ satisfy the following:
they are either both rational integers or both integers in
$\Q(\sqrt{5})$; while $a_\wp = b_\wp$ in the first case, they are
conjugated to each other in the second. In particular, this
implies that the case $L=K$ is impossible, since else we would
have $\Q(\{a_\wp\}) = \Q$.

Suppose then that $[L:K]=4$. This would imply that the
representations $\sigma_\ell$ would be absolutely irreducible.
Indeed, this is a consequence of Faltings' proof of Tate's
conjecture since by theorem \ref{L} we know that in this case
$\End_K(A)= \Z$.

Now, since $\Gal (L/K)=D_2$, we know that there are three
intermediate fields $K\varsubsetneq E_1, E_2, E_3\varsubsetneq L$
satisfying $E_1 \cdot E_2 =
 E_2 \cdot E_3 = E_3 \cdot E_1 = L$. For each of them, we know that
$\End_{E_i}(A)$ is an order $\cR _i$ in a quadratic field $W_i$.

We wish to apply the results of the previous section to the
extensions $L/E_i$ and we first need to explain how these results
generalize to the case of a non real field $W_i$. The only changes
concern the determinant of the two-dimensional irreducible
components $\rho^i_\la$ and ${\rho^i_\la}^{\gamma_i}$ of the
representations  of $\Gal(\qbar / E_i)$: by the Riemann
hypothesis, it easily follows that in the standard factorization
(5.1) of the characteristic polynomials it must hold that $a_\wp ,
b_\wp \in \mathbb{R}$. Thus, if $W_i$ is not real, we  have that
$\tr (\rho^i_\la (\Fr \; \wp)) = a_\wp \in W_i \smallsetminus \Q
\Rightarrow
 a_\wp \not\in \mathbb{R}$ for almost every prime
$\wp$ in $E_i$ inert in $L / E_i$ and hence
 $\det(\rho^i_\la (\Fr \; \wp)) \neq p $. This and the
description of $\rho^i_\la|_{I_\ell}$ forces the determinant of
$\rho^i_\la$ to be equal to $\phi_i \cdot \chi$ for some
non-trivial finite order character $\phi_i$ unramified outside
$N$. Finally, since the determinant is $\chi $ when restricted to
the subgroup $G_L$ of index two of $G_K$, we conclude that
$\phi_i$ is precisely the quadratic character corresponding to
$\Gal(L/E_i)$. Formula (4.5) is thus verified (with the obvious
change of notation: $\psi$ becomes $\phi_i$).

Let us call $T_i = \Gal(\overline{\Q}/E_i)$, $i=1,2,3$. We
conclude that $\sigma_\ell$ is absolutely irreducible and we know
that it contains the reducible groups $\sigma_\ell|_{T_i} =
\rho^i_\la \oplus {\rho^i_\la}^{\gamma_i}$ as normal subgroups of
index two, where $\langle \gamma_i \rangle = \Gal(W_i/\Q)$.

We also know that the extensions $E_i/K$ only ramify at the primes
in $N$ so there are finitely many options for them and thus also
for $L$. Let $H=\Gal(\qbar /L)$. We know that $\sigma_\ell|_H =
(\rho^i_\la \oplus {\rho^i_\la}^{\gamma_i})|_H = \rho_\ell \oplus
\rho_\ell$ for $i=1,2$ and $3$ and $\rho_\ell$ with values in
$\GL(2, \Z_\ell)$ for almost every $\ell$. Let $\wp$ be a prime in
$K$. If $\wp$ is totally decomposed in $L/K$, then $\Fr \; \wp \in
H$ and $\tr (\sigma_\ell (\Fr \; \wp)) = 2 \cdot \tr (\rho_\ell
(\Fr \; \wp)) = 2 a_\wp $ with $a_\wp \in \Z$.

On the other hand, let $\wp$ be a prime in $K$ not totally
decomposed in $L/K$. Despite this fact, there exists $i \in
\{1,2,3\}$ such that $\wp$ decomposes in $E_i/K$, but is inert in
$L/E_i$. Then $\Fr \; \wp \in T_i \smallsetminus H$ and by
applying formula (4.5) we obtain that $\tr (\sigma_\ell (\Fr \;
\wp)) = \tr (\rho^i_\la (\Fr \; \wp)) + \tr
({\rho^i_\la}^{\gamma_i}(\Fr \; \wp)) = a_\wp + a_\wp^{\gamma_i} =
a_\wp + \phi_i(\Fr \; \wp)  a_\wp = 0$.

\begin{lemma}
\label{teo:irreducible}

Let $A/K$ be an abelian surface with quaternionic multiplication
and let $L/K$ be the minimal field of definition of the
endomorphisms of $A$. Let $N$ be the product of the primes of bad
reduction of $A$ over $K$. Then, if $\Gal (L/K)=D_2$, there exist
two different quadratic extensions $E_1$ and $E_2$ of $K$, both
unramified outside $N$, such that

\begin{itemize}
\item
$L$ is the compositum of $E_1$ and $E_2$.

\item
For every prime $\wp \nmid N$ of $K$ totally decomposed in $L/K$,
the characteristic polynomial of $\sigma_\ell (\Fr \; \wp)$, when
factorized as in (5.1), verifies $a_\wp = b_\wp \in \Z$.
\end{itemize}

On the other hand, if $\wp$ does not totally decompose in $L/K$,
then

$$
\tr (\sigma_\ell (\Fr \; \wp))= 0.
$$
\end{lemma}

There is a finite number of possibilities for $E_1$ and $E_2$: in
the example considered, these two fields must be two (different)
extensions of $K =\Q(\sqrt{-3})$ unramified outside $N= 2 \cdot
\tres \cdot 5 \cdot 7 \cdot 11$. Computations show that, for any
choice of such a pair of quadratic extensions, there is a prime
$\wp$ of $K$ not totally decomposed in the compositum field
contradicting the trace  $0$ condition of the lemma above.

We recall that, in order to simplify  computations, we only
computed the characteristic polynomials for primes $\wp$ of $K$
that have residue class degree $1$. We have performed these computations
for all such $\wp$ with residue characteristic $p \leq 193$.
 Therefore, in virtue of lemma
\ref{teo:irreducible}, we conclude that $L/K$ is not a quartic
extension in our example. Having eliminated two of the three cases
of theorem \ref{L}, {\bfseries I}, we conclude that the Jacobian
variety $J(C)$ of Hashimoto- Murabayashi's curve $C$ has
quaternionic multiplication over a quadratic extension $L$ of $K =
\Q(\sqrt{-3})$ and that $\End _K(J(C))$ is the real quadratic
field $\Q(\sqrt{5})$.

The quadratic extension $L/K$ is unramified outside $N= 2 \cdot
\tres \cdot 5 \cdot 7 \cdot 11$ and formula (4.5) tells us that a
non zero trace $a_\wp$ is in $\Z$ if and only if the prime $\wp$
decomposes in $L/K$. Thus, considering all possible quadratic
extensions of $K$ unramified outside $N$ and applying formula
(4.5) to the traces computed, we see that the only extension that
matches is $L= \Q( \sqrt{-3}, \sqrt{-11})$.

\subsection{Explicit determination of the images of the Galois representations}

We now wish to compute the finite set of (possibly) exceptional
primes of the Galois representations on the Tate modules of the
Jacobian variety of Hashimoto-Murabayashi's curve $C$. By theorem
\ref{teo:C2}, we are placed under the assumptions of theorem
\ref{teo:ribetimages} in section 4.

The hard part of the task is determining the primes such that the
residual representation fails to be irreducible or does not have a
multiple of $\ell$ order. These primes must fall in one of the
following cases:

\begin{enumerate}

\item
$\rer$  reducible

\item
$\mathbb{P}(\rer)$  cyclic

\item
$\mathbb{P}(\rer)$ dihedral

\item
$\mathbb{P}(\rer)$ small exceptional
\end{enumerate}

\subsubsection{Reducible primes}

We begin with the determination of those primes falling in cases
1) and 2), i.e., primes such that $\rer$ is reducible over
$\overline{\F}_\la$. We will call them {\em reducible primes}.
Applying Raynaud's result, we see that if $\la \nmid N$ is a
reducible prime, we are in one of the two following situations:

\begin{itemize}
\item[(a)]
$\rer \simeq
\begin{pmatrix}
\epsilon & * \\
0 & \epsilon^{-1} \chi
\end{pmatrix}
$
\item[(b)]
$\rer \simeq  \begin{pmatrix} \epsilon \psi_2 & * \\ 0 &
\epsilon^{-1} \psi_2^{\ell}
\end{pmatrix}
$
\end{itemize}
where $\epsilon$ is, in both cases, a character unramified outside
$N$,  $\chi$ is the $\mod \; \ell$ cyclotomic character and
$\psi_2$ a fundamental character of level $2$.

In order to control the character $\epsilon $ we can use the bound
for conductors of abelian varieties given in \cite{BrKr}. Since
$A$ is an abelian surface defined over $\Q(\sqrt{-3})$, we obtain
that $\cond (A) \mid 2^{20} \cdot \tres^{16} \cdot 5^9 \cdot 7^4
\cdot 11^4$. This is a bound for the conductor of $\re \oplus
\re^{\gamma}$ and thus we can assume that, for any prime $\la$ in
$\Q(\sqrt{-3})$ it holds that $\cond (\re) \mid 2^{10} \cdot
\tres^8 \cdot 5^5 \cdot 7^2 \cdot 11^2$. If $\la$ is a reducible
prime as in case (a) or (b) above, the character $\epsilon$  must
therefore verify: $\cond (\epsilon) \mid 2^5 \cdot \tres^4 \cdot
5^2 \cdot 7 \cdot 11$. Let us first treat

{\em Case (a)}: Equating traces we obtain

$$
( a_\wp \; \mod \; \lambda) = \epsilon(\Fr \; \wp) +
\epsilon^{-1}(\Fr \; \wp) p
$$
as elements in $\overline{\F}_\la$, for every $\wp \nmid \ell N$.
Let $K:= \Q(\sqrt{-3})$ and $\cR$ its ring of integers. We will
apply class field theory over $K$ to compute the reducible primes.
We shall use repeatedly the fact that $K$ has class number $1$ and
that the only units in $K$  are the sixth roots of unity.

Observe that the image of $\epsilon$ is contained in
$\overline{\F}_\la^*$, so this character of the Galois group $G$
of $K$ corresponds to a cyclic extension of $K$ unramified outside
$N$ with conductor dividing $c:= 2^5 \cdot {\tres}^4 \cdot 5^2
\cdot 7 \cdot 11$. If we call

$$
P(c) = \{ \wp \; : \; \mbox{there exists} \; \pi \in \cR \;
\mbox{with} \;  \wp = (\pi) \; \mbox{and} \; \pi \equiv 1 \pmod{c}
\}
$$
and let $F$ denote the ray class field of $K$ of conductor $c$,
then $F$ is characterized by the fact that $F/K$ is abelian and
the set $P(c)$ is exactly the set of prime ideals of $K$ that
decompose totally in $F$. The cyclic extension $F'$ of $K$
corresponding to $\epsilon$ is of course contained in $F$. Thus,
given a prime $\ss = (t)$ in $K$ verifying $ t^f \equiv 1
\pmod{c}$, we have $(\ss^f, F/K)= (\ss , F/K)^f = 1 \in \Gal(F/K)$
and from $K \subseteq F' \subseteq F$ we obtain: $\epsilon(\Fr \;
\ss)^f =1$. From this and the assumption that the characteristic
polynomial of $\rer(\Fr \; \ss)$ admits $\epsilon(\Fr \; \ss)$ as
a root, we obtain the equation for the resultant:

$$
\Res_q := \Res(x^2 - a_{\ss} x + q , x^f -1) \equiv 0 \pmod{\la}
$$
for $\ss \nmid \ell N$, $q$ the rational prime below $\ss$ and
$\ss = (t)$ for some $t$ with $t^f \equiv 1 \pmod{c}$.

In the example, we apply this equation with $q= 31, 43$ and $61$,
$\ss \mid q$ generated by $ t= 2+3 \sqrt{-3}, 4 + 3 \sqrt{-3}$ and
$ 7 + 2 \sqrt{-3}$ (respectively) having all them order $\mod \;
c$ equal to $f= 240$. The values of the traces are $a_{\ss} = -4 ,
4 \sqrt{5} $ and $ 4 \sqrt{5}$ (respect.). Having computed
$\Res_q$ for these three values of $q$, we see that for every
prime $\la$ in $\Q(\sqrt{5})$ with $\la \mid \ell > 11$ one of  them verify
$\la \nmid \Res_q, \ell \neq q$. Thus we conclude that $\rer$ is
not reducible as in case (a) for any $\ell > 11$.

{\em Case (b)}: The analysis made in case (a) tells us how to
control character $\epsilon$, now it remains to say a few words
about the ``fundamental character" $\psi_2$: in fact, we are
abusing notation since we are denoting $\psi_2$ a character of
$G_K$ unramified outside $\ell$
 whose restriction to $I_\ell$ agrees with a level $2$ fundamental character.
 We can identify these two characters because $K$ has class number $1$.
  Let $\ss=(t)$ be a prime in $K$,
then we know that $\psi_2(\Fr \; \ss) \equiv \zeta t'
\pmod{\lambda}$ where $\zeta$ is a unit in $K$ ($\zeta^6 = 1$) and
$t' = t$ or $t^\alpha$, $\alpha$ the order two element in
$\Gal(K/\Q)$. Observe that for case (b) to hold it should be
$\ell$ inert in $K/\Q$, and so $t^\alpha \equiv t^\ell
\pmod{\lambda}$. Take as in the previous discussion of case (a) a
prime $\ss=(t)$ of $K$ with $t^f \equiv 1 \pmod{c}$. Increase $f$
if necessary so that $6 \mid f$. Then, using the assumption that
both $\epsilon(\Fr \; \ss ) \psi_2(\Fr \; \ss)$ and
$\epsilon^{-1}(\Fr \; \ss) \psi_2^\ell(\Fr \; \ss)$ are roots of
the characteristic polynomial of $\rer(\Fr \; \ss)$ we conclude
that the following equation is satisfied

$$
\Res (x^2 - a_{\ss} x + q , x^f - t^f) \equiv 0 \pmod{\lambda}
$$
for $\ss \nmid \ell N$, $q$ the rational prime below $\ss$ and
$\ss = (t)$ for some $t$ with $t^f \equiv 1 \pmod{c}$, $6 \mid f$.
We apply this equation in the example with $q=43, 61$ and $193$.
For $q=193$ we take $\ss \mid q$ generated by $ t= 1+8 \sqrt{-3}$,
this generator has order $\mod \; c$ equal to $f=120$ and the
corresponding trace is $a_{\ss} = 6 \sqrt{5}$. From these
computations it follows that if $\la \mid \ell > 11$ and $\ell
\neq 89$ the residual representation $\rer$ does not fall in case
(b).

It only remains to say a word about $\ell= 89$. In order to prove
that this is not reducible prime, we use the following fact: let
$H= \Gal(\qbar /L)$ be the absolute Galois group of the field $L$
of definition of the quaternionic endomorphisms. We know that, for
any $\lambda \nmid N \cdot D$, the image of the restriction
$\rer|_H$ lies in $\GL_2(\F_\ell)$. Combined with the assumption
in case (b), this shows that (the semisimplification of) $\rer|_H$
is contained in a non-split Cartan subgroup of $\GL_2(\F_\ell)$ (
$\mathbb{P}(\rer|_H)$ is cyclic). Thus, the image of $\rer|_H$ can
contain no matrix whose characteristic polynomial is reducible
over $\F_\ell$ with two different $\mod \; \ell$ eigenvalues. We
have computed a few characteristic polynomials for primes $\wp$
such that $\Fr \; \wp \in H$ (recall that we have shown that this
is the case if $a_\wp \in \Z$, $a_\wp \neq 0$) and we found that,
for $p=157$, the corresponding characteristic polynomial is
$x^2+4x+157$ and that it reduces $\mod \; 89$ with two different
eigenvalues. This shows that  $89$ is not a reducible prime as in
case (b).

\subsubsection{Dihedral and small exceptional primes}

The determination of dihedral primes is carried out by using the
technique applied in \cite{Se}, \cite{Ri5} and \cite{DiVi} via the
description of the restriction to $I_\ell$ provided by Raynaud's
theorem.

If $\mathbb{P}(\rer)$ is dihedral, then there exists a Cartan
subgroup $C$ such that the image $\overline{\V}$ of $\rer$ is
contained in the normalizer $\mathcal{N}$ of $C$ but not in $C$ itself.
Composing $\rer$ with the quotient $\mathcal{N}/C \simeq C_2 \simeq \{ \pm 1
\}$, we obtain a quadratic character $\phi$ of the Galois group
$G=\Gal (\qbar /K)$ corresponding to a quadratic extension
$E_\ell$ of $K$ unramified outside $\ell N$. Furthermore, the
description of $\rer|_{I_\ell}$ for $\ell \nmid N$ shows that, if
$\ell>3$, it must be contained in $C$, so $E_\ell /K$ does not
ramify at $\ell$. The traces of elements in $\mathcal{N} \smallsetminus C$
all vanish and the value $\phi(\Fr \; \wp)$, where $\phi$ only
ramifies at the prime divisors of $N$, determines whether
$\rer(\Fr \; \wp)$ falls in $C$ or not. Thus, for every prime
$\wp$ in $K$ such that $\wp$ is inert in $E_\ell /K$, we have that
$\phi(\Fr \; \wp) = -1$ and hence $a_\wp \equiv 0 \pmod{\la}$.

Therefore, the algorithm to compute all dihedral primes is the
following:

\begin{itemize}

\item
List all quadratic extensions $E$ of $K$ unramified outside
 $\{2,\tres,5,7, 11 \}$.

\item
For each of these extensions, find several primes $\wp$ in $K$
such that

$$
\wp \; \mbox{inert in} \; E/K  \quad \mbox{and} \quad  a_\wp \neq
0 \quad \; \eqno(\dag )
$$
\end{itemize}

If $\la \mid \ell$ is a dihedral prime, we then should have that,
for some quadratic extension $E/K$ as above and all primes $\wp
\nmid \ell $ verifying $\dag $,

$$
\la \mid a_\wp .
$$

We have computed the traces $a_\wp$ for every prime $\wp$ of $K$
with residue class degree $1$ and $\Norm(\wp) \leq 193$ and
applied the above algorithm to all quadratic extensions of
$\Q(\sqrt{-3})$ unramified outside $\{ 2, \tres , 5, 7, 11 \}$ and
we found no dihedral primes $\lambda \mid \ell >11$. Finally, in
order to eliminate the possibility of primes with small
exceptional image, we use the trick applied in \cite{Ri5},
\cite{DiVi}: the description of $\mathbb{P}(\rer|_{I_\ell})$ shows
that this group contains a cyclic subgroup of order $\ell \pm 1$.
Hence, if $\ell > 5$, $\ell \nmid N$, the image of
$\mathbb{P}(\rer)$ can not be isomorphic to neither $A_4 , S_4$
nor $A_5$. We have thus proved the following

\begin{theorem}

Let $C/K$ be a smooth projective model of the curve

$$
Y^2 = \frac{1}{48}X(9075X^4+3025(3+2 \sqrt{-3}) X^3 -6875 X^2 +220
(-3+ 2 \sqrt{-3}) X + 48)
$$
over $K=\Q (\sqrt{-3})$ and let $A/K$ be its Jacobian variety. Let
$\{ \re \}$ be the two-dimensional Galois representations on the
Tate modules of $A$. Then, for every prime $\ell >11$, $\la \mid
\ell$, the residual representation $\rer$ is absolutely
irreducible and the order of its image is a multiple of $\ell$.

\end{theorem}

In order to apply theorem \ref{teo:ribetimages}, we just observe
that, for primes $\wp$ above $13$ and $19$, we have $a_\wp = 2
\sqrt{5}$, whereas for primes $\ss$ above $31$ and $37$, we have
$a_{\ss} = -4$, $4$, respectively.

\begin{theorem} For every prime
$\ell >11$, $\la \mid \ell$, the image of the Galois
representation $\re$ is the subgroup of $\GL_2(\mathbb{Z}[\sqrt{5}]_\la)$ generated
by $\GL_2(\Z_\ell)$ and the diagonal matrix $\begin{pmatrix}
  \sqrt{5} & 0 \\
  0 & \sqrt{5}^{-1}
\end{pmatrix} $. In particular, if $\ell \equiv \pm 1 \pmod{5}, \; \ell \neq 11$,
the groups $\GL_2(\Z_\ell)$  and $\GL_2(\F_\ell)$ are realized as
Galois groups over  $\Q(\sqrt{-3})$ and the corresponding
extension is unramified outside $2310 \ell$. Furthermore, for
every prime $\ell >11$, the group $\PGL_2(\F_\ell)$ is realized as
a Galois group over $\Q(\sqrt{-3})$, again through an extension
unramified outside $2310 \ell$.

\end{theorem}

\end{document}